\documentclass[%
 aip,
 cha,
 amsmath,amssymb,
 reprint,%
]{revtex4-1}
\usepackage{lmodern}

\usepackage{dcolumn}
\usepackage{bm}

\usepackage[utf8]{inputenc}
\usepackage[T1]{fontenc}
\usepackage{mathptmx}
\usepackage{etoolbox}
\usepackage{amssymb, amsmath, amsthm, amsfonts}

\usepackage[cm]{fullpage}
\usepackage[english]{babel}
\usepackage[T1]{fontenc}
\usepackage[pdftex]{graphicx}
\usepackage{booktabs}
\usepackage{xcolor}
\usepackage{xfrac}

\usepackage{mathtools}
\usepackage{epstopdf}
\usepackage{dsfont}
\usepackage{mathrsfs}
\usepackage{physics}
\usepackage{bm}
\usepackage{diagbox}

\usepackage{hyperref}
\hypersetup{
    colorlinks = false
}

\newcommand{\R}{\mathbb{R}}

\newcommand{\E}[1]{\mathbb{E}\left( #1 \right)}
\newcommand{\PRo}[1]{\mathbb{P}\left( #1 \right)}

\newcommand{\lap}[2][t]{\mathscr{L}_{#1}\!\left\{ #2 \right\}\!}
\newcommand{\four}[2][t]{\mathscr{F}_{#1}\!\left\{ #2 \right\}\!}

\DeclareMathOperator*{\argmax}{arg\,max}

\makeatletter
\def\@email#1#2{%
 \endgroup
 \patchcmd{\titleblock@produce}
  {\frontmatter@RRAPformat}
  {\frontmatter@RRAPformat{\produce@RRAP{*#1\href{mailto:#2}{#2}}}\frontmatter@RRAPformat}
  {}{}
}%
\makeatother
\begin{document}

\preprint{AIP/123-QED}

\title{Arcsine laws for Brownian motion with Poissonian resetting}
\author{K. Taźbierski}
\email{kacper.tazbierski@pwr.edu.pl}
\author{M. Magdziarz}%
\affiliation{Faculty of Pure and Applied Mathematics, Wroc{\l}aw University of Science and Technology,\\ Wyb. Wyspia{\'n}skiego 27, 50-370 Wroc{\l}aw, Poland}%

\date{\today}

\begin{abstract}
We analyze the equivalents of the celebrated arcsine laws for Brownian motion undergoing Poissonian resetting. We obtain closed-form formulae for the probability density functions of the corresponding random variables in the cases of the first and second arcsine law. Furthermore, we obtain numerical results for the third law.
\end{abstract}

\maketitle

\begin{quotation}
The celebrated arcsine laws describe surprising statistical properties of Brownian motion, capturing the distribution of specific random variables associated with its temporal evolution. In this work, we explore analogous phenomena for Brownian motion under Poissonian resetting, a process that introduces interruptions at random times, fundamentally altering its dynamics. We derive closed-form expressions for the probability density functions corresponding to the first and second arcsine laws, highlighting their distinct behavior under resetting. Additionally, we present numerical insights into the third arcsine law, providing a comprehensive perspective on this generalization of classical results.
\end{quotation}
\section{Introduction}
One of the basic models in the theory of stochastic processes is Brownian motion, whose applications spread over vast areas, both in the natural and social sciences: physics, finance, and biology. This process, also known as the Wiener process, describes the random diffusive motion of a particle floating in a fluid, characterized by paths that are continuous, nowhere differentiable with linear in time mean squared displacement. These properties make it a prominent mathematical model of various phenomena and a central object of study in stochastic processes \cite{Einstein1905,Wiener1923,Morters2010}.

Among many properties of Brownian motion the arcsine laws are particularly surprising and elegant \cite{Yen2013}. They give information about distribution of the time-related quantities associated with the process. These laws appear in several contexts, including random walks \cite{Alsmeyer2019}. More precisely, the arcsine laws describe the distribution of the proportion of time that a Brownian particle spends above a particular level, as well as the distribution of the time at which the process reaches its maximum and the last time the process hits zero on the time interval $[0,1]$.

The first arcsine law describes the distribution of the proportion of time that one-dimensional Brownian motion spends above zero on the interval $[0,1]$ (see Section \ref{Sec2} for the detailed formulation). This random quantity is arcsine distributed. This surprising result implies that the Brownian particle is more likely to spend most of its time either above or below zero, rather than to spend half the time each on the positive and negative half-line \cite{Kotani2006,Yor1992,Spitzer1956}.

The second arcsine law describes the last moment before \(t=1\) in which Brownian particle hits zero \cite{Yen2013}. This moment, again, is arcsine distributed, which implies that the event is more likely to have taken place closer to the beginning or the end of $[0,1]$ than in the middle of this time interval \cite{Knight1963,Borodin2002,Gallardo2021}. 

The third arcsine law describes the distribution of the time at which Brownian particle reaches its maximum on the interval $[0,1]$ \cite{Yen2013}. This random time is arcsine distributed as well, meaning that the maximum is probably reached either very early or very late on the time interval \([0, 1]\). The result is counterintuitive because a priori one would expect the maximum to be attained near the mid-time of the interval \cite{Knight1963,Karatzas1991,Ding2019}.

All three laws put together shed light on the unique and often surprising behavior of Brownian motion, and hence are of particular importance for both theoretical and applied probability \cite{Revuz1999,Billingsley1995,Breen2022}.

In recent years, the idea of resetting in the context of stochastic dynamics has become a major field of study in both mathematics and statistical physics \cite{Evans2011, Evans2020}. Stochastic resetting means that at certain random times the diffusing particle is taken back to some fixed location, restarting the dynamics anew. This idea has attracted great attention since resetting can significantly change the statistical properties of the underlying process, leading to new and interesting behavior of complex systems, see \cite{Evans2011, Evans2020} and the references therein.

The stochastic resetting mechanism has been extensively used in the context of search processes and reaction kinetics, where it is able to optimize search efficiency by preventing excessive excursions away from the target \cite{Evans2011, Pal2015}. The interplay of Brownian motion with stochastic resetting has led to the finding of new classes of non-equilibrium stationary states \cite{Reuveni2016, Ray2020}. Moreover, many recent works have concerned first-passage times for processes with resetting. It turns out that particles under resetting are able to find the target much faster than expected under many conditions \cite{Pal2017,Besga2020}. These results have immediately found application in various contexts and disciplines, such as: ecology, animal foraging strategies or computer science, where resetting could significantly improve searching procedures with algorithms \cite{Pal2015, Kusmierz2014}. One can find more information on the theory and applications in an introductory article by Kundu and Reuveni \cite{Kundu2024}. A popular type of stochastic resetting, called sharp resetting, that involves resetting a stochastic process (e.g., a particle's motion) at fixed, deterministic time intervals, quickly attracted attention due to its simplicity and power as a modeling tool. It was also extensively studied in \cite{Eli2020, Eli2021, Eli2021_2, Eli2023, Eli2023_2}. Unlike general stochastic resetting, which uses probabilistic intervals for reset events, sharp resetting enforces strict periodicity in the timing of resets. This approach is used to optimize search strategies, improve reaction rates and control random processes in systems that benefit from regular interruptions. Many different approaches to resetting were studied, such as the so-called stochastic ratcheting \cite{Ghosh2023}. There, the resetting mechanism introduces spatial symmetry breaking via a certain external potential and introduces many resetting points. Discontinuous behavior of the process is also eliminated. This method allows to model movements of bacteria or micro-robots which turn their engines on and off.

In this paper, we will study the analogues of arcsine laws for Brownian motion under stochastic resetting. We will derive analytical results for the distribution of the proportion of time that the Brownian particle under resetting spends above zero. We will also find the distribution of the last time the particle hits zero.
\section{First (L\'evy's) arcsine law}
\label{Sec2}
Let $T$ be a total time that the Brownian motion (Wiener process) $W(t)$
spends on the positive half-line during the time interval $[0,1]$, i.e.
\begin{equation}
    T=\lambda\left(\{t\in[0,1]: W(t)\geq 0  \}\right),
\end{equation}
where $\lambda$ is the Lebesgue measure. The Wiener process may be characterized through a set of properties \cite{Karatzas1991} 
\begin{enumerate}
    \item $W(0)=0$,
    \item $W$ has independent and stationary increments,
    \item $W(t)$ is normally distributed with mean $0$ and variance $t$,
    \item $W$ is continuous almost surely.
\end{enumerate}
The first arcsine law states that $T$ is arcsine distributed, i.e.
\begin{equation}\label{eq:Tr}
    \PRo{T\leq t}=\frac{2}{\pi}\arcsin{(\sqrt{t})}
\end{equation}
for $t\in[0,1]$.

In what follows, we will make extensive use of the important scaling property
\begin{gather}
    \nonumber\PRo{\lambda\left(\{\eta\in[0,\theta] : W(\eta)\geq 0  \}\right)\leq t}\\ \label{eq:Tscaling}=\PRo{\lambda\left(\{\eta\in[0,1] : W(\eta)\geq 0  \}\right)\leq {t}/{\theta}}.
\end{gather}

Now let $X(t)$ be a Wiener process starting at 0 undergoing Poissonian resetting with intensity $r$ to point $x_r=0$ \cite{Evans2020}. Recall that the Poisson process with intensity $r>0$  is a counting process $N(t)$ fulfilling \cite{feller1957}
\begin{enumerate}
    \item $N(0)=0$,
    \item $N$ has independent and stationary increments,
    \item $N(t)$ is Poisson distributed with mean $rt$
\end{enumerate}
We also know that the random vector $\bm{\tau}'=(\tau'_1,\ldots,\tau'_k)$, where $\tau'_i=\left(\tau_i |N(1)=k\right)$ and $\tau_i$ is the time of the $i$th jump of $N(t)$, has the PDF \cite{feller1957}
\begin{equation}\label{eq:poissonuniform}
    p_{\bm{\tau}'}(t_1,\ldots,t_k)=k!\mathds{1}(0< t_1<\dots< t_k< 1),
\end{equation}
where $\mathds{1}$ is the indicator function. Here, the notation with "$|$" means that we condition on the event $N(1)=k$. This implies that the $i$th jump time of the Poisson process on the interval $[0,t]$ conditioned on $N(t)=k$ is distributed as the $i$th order statistic of a vector of $k$ uniformly $\mathcal{U}(0,t)$ distributed random variables \cite{ross2007}.

Denote by $T_r$ a total time that $X(t)$ spends on the positive half-line during the time interval $[0,1]$, i.e.
\begin{equation}
    T_r=\lambda\left(\{t\in[0,1]: X(t)\geq 0  \}\right).
\end{equation}
In what follows we will derive the formula for the probability density function (PDF) $p_{T_r}(t)$ of $T_r$.

\begin{figure*}
    \centering
    \includegraphics[width=\linewidth]{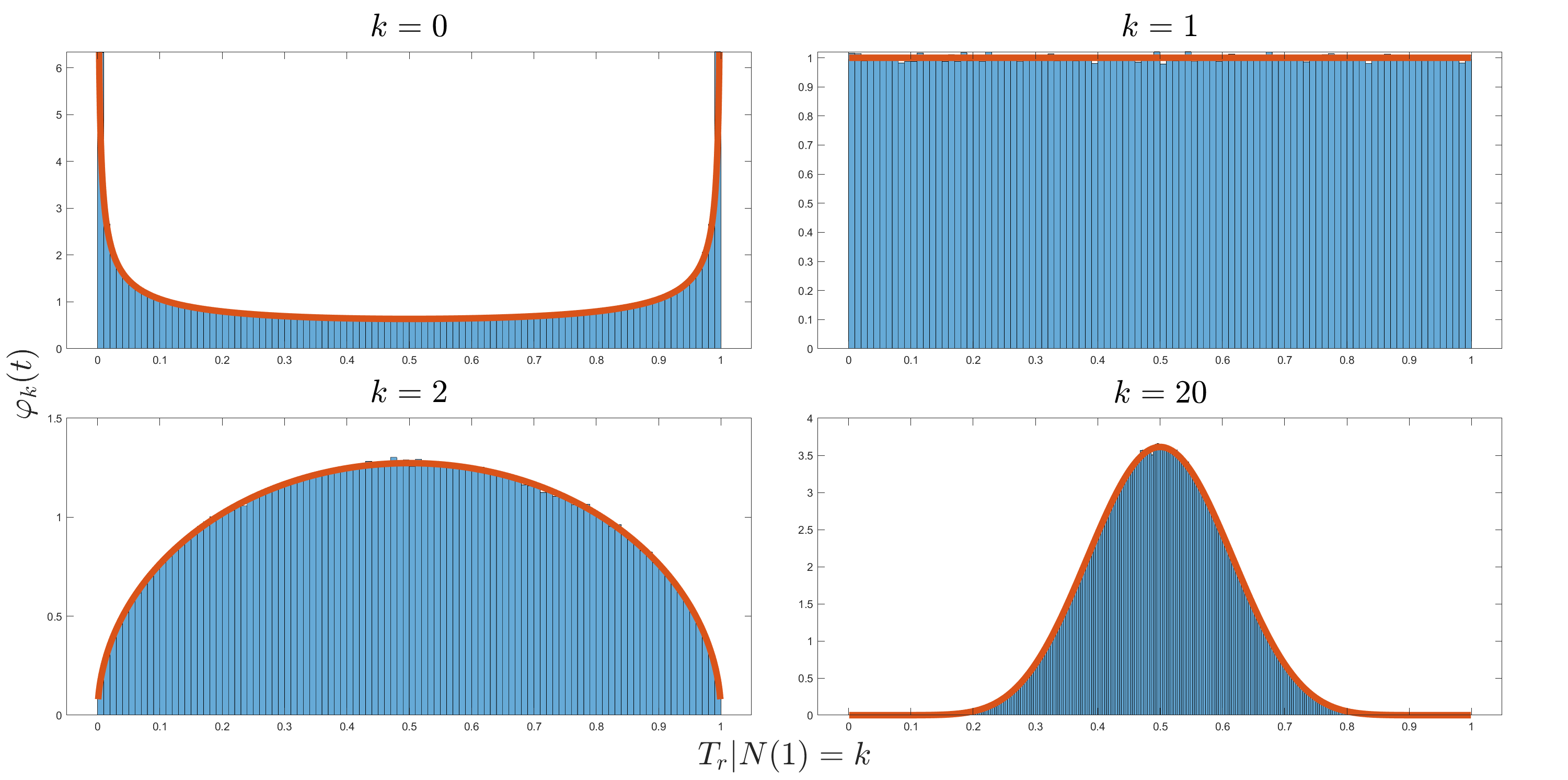}
    \caption{Monte-Carlo estimated distribution of $T_r|N(t)=k$ \eqref{eq:Tr}, i.e. the time the Wiener process under resetting spends above the $X$-axis when $k$ resets occurred,  with sample size $10^6$ compared to the theoretical result \eqref{eq:TKFinal}. The theoretical PDF seems to  perfectly fit the simulations. With increasing number of resets $k$ we can observe the convergence of the distribution to a point-mass concentrated at $\sfrac{1}{2}$.}
    \label{fig:Tk}
\end{figure*}

\begin{figure*}
    \centering
    \includegraphics[width=\linewidth]{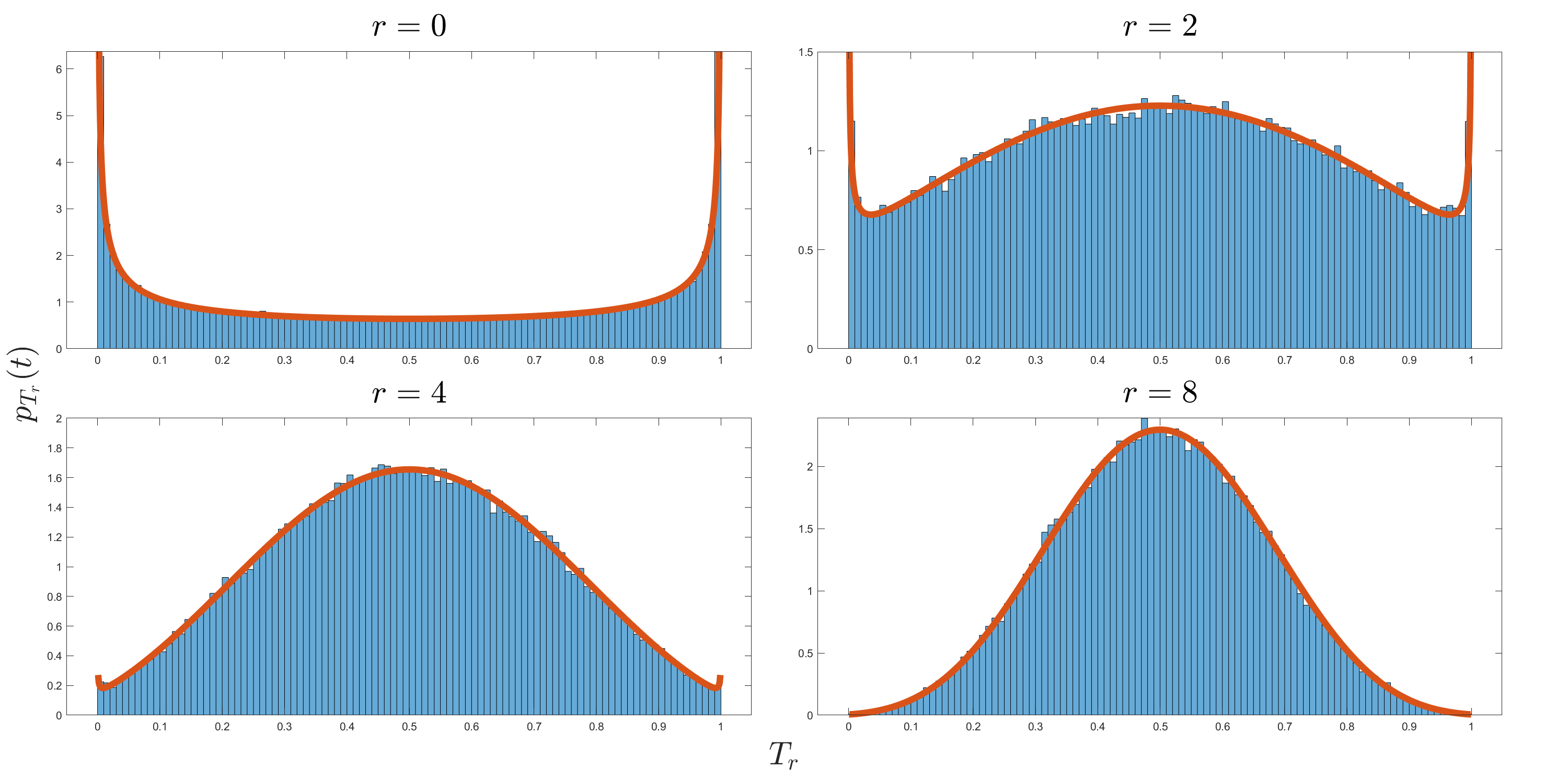}
    \caption{Monte-Carlo estimated distribution of $T_r$ \eqref{eq:Tr}, i.e. the time the Wiener process under resetting with rate $r$ spends above the $X$-axis,  using Wiener process trajectories with sample size $10^5$ and timestep $10^{-4}$. In red we see the theoretical result \eqref{eq:TrResult}. The quality of the fit visually seems perfect. In accordance with the case presented on figure \ref{fig:Tk}, when the resetting rate increases, the distribution approaches a point-mass at $\sfrac{1}{2}$.}
    \label{fig:Tr}
\end{figure*}

On the interval $[0,1]$ we have $N(1)\sim Poiss(r)$ resetting events, with resetting moments uniformly distributed on $[0,1]$ \cite{ross2007}. Thus, we have from the total probability formula
\begin{align}\label{eq:Ttotal}
    \nonumber p_{T_r}(t)\dd t&=\PRo{T_r\in(t,t+\dd t)}\\
    &=\sum_{k=0}^\infty \PRo{T_r\in(t,t+\dd t)|N(1)=k}\PRo{N(1)=k}.
\end{align}
Now let us look at the random variable $T_r|N(1)=k$. The interval $[0,1]$ is split into $k+1$ parts. In each subinterval there is a single independent realization of the Wiener process and the times spent in the positive half-line corresponding to each of them sum up to give us the result for $T_r$. Using the scaling property \eqref{eq:Tscaling} we now get
\begin{equation}\label{eq:TkDef}
    (T_r|N(1)=k)\stackrel{D}{=}\sum_{i=0}^{k}(\tau'_{i+1}-\tau'_i)T^i,
\end{equation}
where $T^i$ is a sequence of iid arcsine-distributed random variables. Here, "$\stackrel{D}{=}$" means equality in distribution. The cases $i=0$ and $i=k+1$ are exempt from the definition and defined as $\tau'_0=0$ and $\tau'_{k+1}=1$. Now we average the term $T_r|N(1)=k$ over the vector $\bm{\tau}'$ using \eqref{eq:poissonuniform}, obtaining
\begin{align}
\nonumber \varphi_k(t)&\equiv p_{T_r|N(1)=k}(t)\\
&=k!\idotsint\limits_{0< t_1<\dots< t_k< 1} p_{\sum_{i=0}^{k}(t_{i+1}-t_i)T^i}(t) \dd t_1\dots\dd t_k.
\end{align}
Now we shift to the Fourier space and get that the Fourier transform $\hat{\varphi}_k(\omega)$ of $\varphi_k(t)$ 
\begin{equation}
    \four[t]{\varphi_k(t)}(\omega)\equiv\int\limits_\R e^{i\omega t}\varphi_k(t)\dd t=\hat{\varphi}_k(\omega),
\end{equation}
equals
\begin{equation}\label{eq:TkCF}
\hat{\varphi}_k(\omega)=k!\idotsint\limits_{0< t_1<\dots< t_k< 1} \prod_{i=0}^k \hat{p}_{(t_{i+1}-t_i)T^i}(\omega) \dd t_1\dots\dd t_k,
\end{equation}
Here we used the independence between $T^i$. Knowing, that $aT$ is generalized-arcsine distributed with support $[0,a]$ we have that
\begin{equation}
    \hat{p}_{aT}(\omega)=e^{\frac{ai\omega}{2}}J_0\left(\frac{a\omega}{2}\right),
\end{equation}
where $J_0$ is the Bessel function of the first kind of order $0$, defined for non-negative real order $\alpha$ as \cite{Abramowitz1964}
\begin{equation}
    J_\alpha(x)=\sum_{k=0}^\infty \frac{(-1)^k}{k!\,\Gamma{\left( 
k+\alpha+1 \right)}}\left( \frac{x}{2} \right)^{2k+\alpha}.
\end{equation}
Plugging this into equation \eqref{eq:TkCF} we get
\begin{equation}\label{eq:T|kCF}
\hat{\varphi}_k(\omega)=k!e^{\frac{i\omega}{2}}\idotsint\limits_{0< t_1<\dots< t_k< 1} \prod_{i=0}^k J_0\left( \frac{\omega(t_{i+1}-t_i)}{2} \right)\dd t_k.
\end{equation}
Changing the variables in the integral with respect to $t_1$ in the following way
\begin{equation}
    u=\frac{\omega t_1}{2}\implies \dd t_1=\frac{2}{\omega}\dd u,
\end{equation}
we get
\begin{widetext}
\begin{equation}
\int\limits_0^{t_2}J_0\left( \frac{\omega t_1}{2} \right)J_0\left( \frac{\omega(t_2-t_1)}{2} \right)\dd t_1=\frac{2}{\omega}\int\limits_0^{\frac{\omega t_2}{2}}J_0\left( u \right)J_0\left( \frac{\omega t_2}{2}-u \right)\dd u=\frac{2}{\omega}J_0^{(2)*}\left(\frac{\omega t_2}{2}\right),
\end{equation}
\end{widetext}
where we denote $f^{(j)*}$ as the $j$-fold convolution of a function with itself. Looking at the multiple integral we see that it consists of a $(k+1)$-fold convolution of $J_0$ functions at point $\omega/2$
\begin{equation}
    \hat{\varphi}_k(\omega)=k!e^{\frac{i\omega}{2}}2^k\omega ^{-k}J_0^{(k+1)*}\left(\sfrac{\omega}{2}\right).    
\end{equation}
Looking at the convolution in Laplace space, we get
\begin{align}
    \nonumber\lap[\omega]{J_0^{(k+1)*}\left(s\frac{\omega}{2}\right)}(s)&=\left(\lap[\omega]{J_0\left(\sfrac{\omega}{2}\right)}(s)\right)^{k+1}\\
    &=\left( s^2+\frac{1}{4} \right)^{-\frac{k+1}{2}},
\end{align}
where the Laplace transform can be checked using power series definition of the Bessel function together with the formula for geometric summation. Inverting the Laplace transform yields \cite{Abramowitz1964}
\begin{equation}\label{eq:JConv}
    J_0^{(k+1)*}\left(\sfrac{\omega}{2}\right)=\frac{\sqrt{\pi}\omega^{\frac{k}{2}}}{\Gamma\left( \frac{k+1}{2}  \right)}J_{\frac{k}{2}}\left( \sfrac{\omega}{2} \right).
\end{equation}
Plugging \eqref{eq:JConv} into \eqref{eq:T|kCF} we now get
\begin{equation}\label{eq:T|kCFinal}
    \hat{\varphi}_k(\omega)=\frac{k!2^k\sqrt{\pi}e^{\frac{i\omega}{2}}}{\Gamma\left( \frac{k+1}{2}  \right)\omega^{\frac{k}{2}}}J_{\frac{k}{2}}\left( \frac{\omega}{2} \right).
\end{equation}
We now may invert this function obtaining \cite{Abramowitz1964}
\begin{align}\label{eq:TKFinal}
    \nonumber\varphi_k(t)&=\frac{k!}{\left(\Gamma\left( \frac{k+1}{2}  \right)\right)^2}\left( \sqrt{t(1-t) }\right)^{k-1}\\
    &=\frac{1}{B\left(\frac{k+1}{2},\frac{k+1}{2}\right)}t^{\frac{k+1}{2}-1}(1-t)^{\frac{k+1}{2}-1},
\end{align}
where $B$ is the Beta function
\begin{equation}
    B(\alpha,\beta)=\int_0^1t^{\alpha-1}(1-t)^{\beta-1}\dd t=\frac{\Gamma{(\alpha)}\Gamma{(\beta)}}{\Gamma{(\alpha+\beta)}}.
\end{equation}
This result has been checked using Monte-Carlo methods as presented in figure \ref{fig:Tk}, yielding the perfect agreement between theory and simulations. Let us notice that this means that $T_r|N_1=k\sim \mathcal{B}\left( \frac{k+1}{2},\frac{k+1}{2} \right)$, i.e. that $T_r$ is conditionally Beta distributed. Now using \eqref{eq:TKFinal} and \eqref{eq:Ttotal} we get
\begin{align}\label{eq:TrResult}
    \nonumber p_{T_r}(t)&= e^{-r}\left( 
t(1-t) \right)^{-\sfrac{1}{2}}\sum_{k=0}^\infty \frac{\left(r\sqrt{t(1-t)}\right)^k}{\left( \Gamma{\left( \frac{k+1}{2} \right)} \right)^2} \\ \nonumber
&=e^{-r}\left( 
t(1-t) \right)^{-\sfrac{1}{2}}\bigg{(} r\sqrt{t(1-t)}J_0\left( 2r\sqrt{t(1-t)} \right) \\
&+ {}_1\bar{F}_2\left(1;\sfrac{1}{2},\sfrac{1}{2};r^2t(1-t)\right)\bigg{)},
\end{align}
where ${}_p\bar{F}_q$ is the regularized generalized hypergeometric function defined as
\begin{equation}
    {}_p\bar{F}_q(\bm{a};\bm{b};z)=\frac{{}_pF_q(\bm{a};\bm{b};z)}{\Gamma{(b_1)}\cdot \ldots\cdot \Gamma{(b_q)}}
\end{equation}
and ${}_pF_q$ is the generalized hypergeometric function
\begin{equation}
    {}_pF_q(\bm{a};\bm{b};z)=\sum_{n=0}^\infty\frac{a_1^{(n)}\cdot\ldots\cdot a_p^{(n)}}{b_1^{(n)}\cdot\ldots\cdot b_q^{(n)}}\frac{z^n}{n!}.
\end{equation}
Also, here $\cdot^{(n)}$ is the Pockhammer symbol, or a rising factorial
\begin{equation}
    x^{(n)}=\prod_{k=0}^{n-1}(x-k).
\end{equation}
This result has again been checked in figure \ref{fig:Tr}, we obtained the perfect agreement between theory and simulations.

As a side note, we present another
approach leading us to exactly the same result by applying a series of naive approximations. Let us return to equation \eqref{eq:TkCF} and check the properties of each term. This random variable is a convex sum of arcsine distributed random variables. The arcsine distribution is a special case of the beta distribution, namely $\mathcal{B}(\sfrac{1}{2},\sfrac{1}{2})$. We also know that the increments of the order statistic of $\mathcal{U}(0,1)$-distributed random variables are beta distributed. The fact that the resetting times can be drawn from the uniform distribution $\mathcal{U}(0,1)$ and then sorted, gives us the first naive approximation - we set $\tau_{i+1}-\tau_i=\frac{1}{k+1}$. Then we have
\begin{equation}\label{eq:Tapr1}
    (T_r|N(1)=k)\stackrel{D}{\approx}\frac{1}{k+1}\sum_{i=0}^{k}T^i.
\end{equation}
Now for the second simplification we will use the method of moments to approximate the sum \eqref{eq:Tapr1} by a beta-distributed random variable $B_k$. The method of moments gives us $B_k\sim\mathcal{B}(\frac{k+1}{2},\frac{k+1}{2})$, thus arriving at the precise analytical result derived previously.

\begin{figure}
    \includegraphics[width=\linewidth]{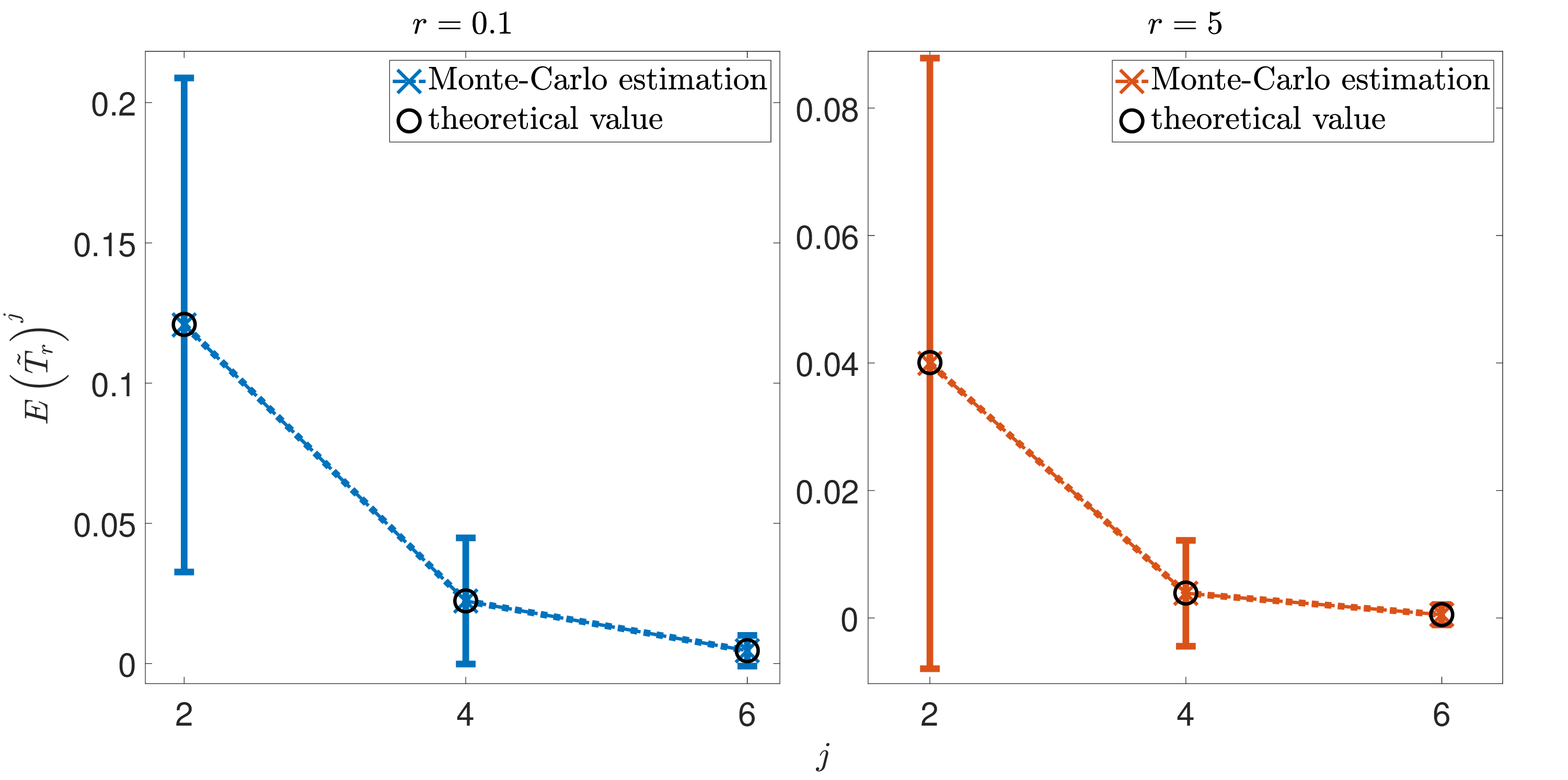}
    \caption{Monte-Carlo estimated $j$th moments of $\tilde{T}_r$ \eqref{eq:Tr} (crosses), their theoretical values \eqref{eq:TrMom} and the standard deviation (whiskers) of the Monte Carlo simulations. We can notice a substantial decrease in the standard deviation due to $\abs{\tilde{T}_r}\leq 1$. }
    \label{fig:Tmom}
\end{figure}

Now, let us find the $j$th central moments of the distribution of $T_r$. Using \eqref{eq:T|kCFinal} with \eqref{eq:Ttotal} we obtain
\begin{equation}\label{eq:TrCF}
    \hat{p}_{T_r}(\omega)=\sqrt{\pi}e^{-r+\frac{i\omega}{2}}\sum_{k=0}^\infty \frac{r^k}{\Gamma\left(\frac{k+1}{2}  \right)\omega^{\frac{k}{2}}}J_{\frac{k}{2}}\left( \sfrac{\omega}{2} \right).
\end{equation}
Using the definition of the Bessel function we can now write \eqref{eq:TrCF} as
\begin{equation}
    \sqrt{\pi}e^{-r+\frac{i\omega}{2}}\sum_{k=0}^\infty \frac{r^k}{\Gamma\left(\frac{k+1}{2}  \right)}\sum_{m=0}^\infty \frac{(-1)^m}{m!\Gamma(m+k/2+1)}\left( \frac{\omega}{2^k} \right)^{2m}.
\end{equation}
Now analyzing a shifted variable $\tilde{T}_r=T_r-\E{T_r}=T_r-\sfrac{1}{2}$ (the mean is a consequence of the characteristic function being purely real), where the exponential coefficient disappears, we can differentiate $\hat{p}_{\tilde{T}_r}(\omega)$ $j$ times to get
\begin{equation}
    \sqrt{\pi}e^{-r}\sum_{k=0}^\infty \frac{r^k}{\Gamma\left(\frac{k+1}{2}  \right)}\sum_{m=j/2}^\infty \frac{(-1)^m(2m)!}{m!\Gamma(m+k/2+1)(2m-j)!}\left( \frac{\omega}{2^k} \right)^{2m}.
\end{equation}
We see that for odd $j$ this sum will yield 0 when $\omega\to0$, so the variable is symmetric around its mean. Analyzing even $j$ we now go with $\omega\to 0$ and obtain
\begin{equation}
    \frac{\dd^j }{\dd \omega^j}\hat{p}_{\tilde{T}_r}(0)=i^j\frac{e^{-r}\Gamma\left(\frac{j+1}{2}  \right)}{2^j}\sum_{k=0}^\infty \frac{\left(\sfrac{r}{2}\right)^k}{\Gamma\left( \frac{k+1}{2} \right)\Gamma\left( \frac{k+j}{2}+1 \right)}.
\end{equation}
Now we can use a well-known formula for the $j-th$ moment of a random variable using characteristic function, to get
\begin{align}\label{eq:TrMom}
    \nonumber\E{(\tilde{T}_r)^j}&=\frac{e^{-r}\Gamma\left(\frac{j+1}{2}  \right)}{2^j}\sum_{k=0}^\infty \frac{\left(\sfrac{r}{2}\right)^k}{\Gamma\left( \frac{k+1}{2} \right)\Gamma\left( \frac{k+j}{2}+1 \right)}\\
    \nonumber &=\frac{e^{-r}\Gamma\left(\frac{j+1}{2}  \right)}{2^j}\Bigg{(}\left(\frac{2}{r} \right)^{\frac{j-1}{2}}I_{\frac{j+1}{2}}(r)\\
    &+{}_1\bar{F}_2(1;\sfrac{1}{2},1+\sfrac{j}{2};(\sfrac{r}{2})^2)\Bigg{)},
\end{align}
where $I_\alpha$ is the modified Bessel function of the first kind of order $\alpha$ defined as \cite{Abramowitz1964}
\begin{equation}
    I_\alpha(x)=i^{-\alpha}J_\alpha(ix).
\end{equation}
In table \ref{tab:table2} we compared the theoretical moments with Monte-Carlo simulations. We calculated the corresponding relative error
\begin{equation}\label{eq:err}
    \varepsilon_j(N,\Delta t)=\abs{\frac{m_j(N,\Delta t)}{\E{\left( \tilde{T}_r \right)^j}}-1}.
\end{equation}
Here, $m_j$ is the estimate of the $j$th central moment with Monte-Carlo sample size $N$ and time step $\Delta t$. The error is of order at most $10^{-3}$ for $N=10^5$, $\Delta t=10^{-5}$ which shows perfect agreement between theory and simulations. Similar comparison is shown in figure \ref{fig:Tmom}. Again we observe a very good fit. 
\begin{table}
\caption{\label{tab:table2}Relative error \eqref{eq:err} of Monte-Carlo approximated moments of $\tilde{T}_r$ using sample size $N=10^5$ with $\Delta t=10^{-5}$ compared to theoretical result \eqref{eq:TrMom}.}
\begin{tabular}{|c|c|c|c|}
\hline
\diagbox[width=6em]{$r$}{moment}& $\E{(\tilde{T}_r)^2}$ & $\E{(\tilde{T}_r)^4}$ &$\E{(\tilde{T}_r)^6}$\\
\hline
$\sfrac{1}{5}$ & $0.52\cdot10^{-3}$& $0.04\cdot10^{-3}$& $1.33\cdot10^{-3}$\\
\hline
2 & $4.53\cdot10^{-3}$ & $3.72\cdot10^{-3}$ & $1.56\cdot10^{-3}$\\
\hline
5 & $0.42\cdot10^{-3}$ & $1.14\cdot10^{-3}$ & $1.75\cdot10^{-3}$\\
\hline
\end{tabular}
\end{table}
\section{Second arcsine law}
\begin{figure}
    \centering
    \includegraphics[width=\linewidth]{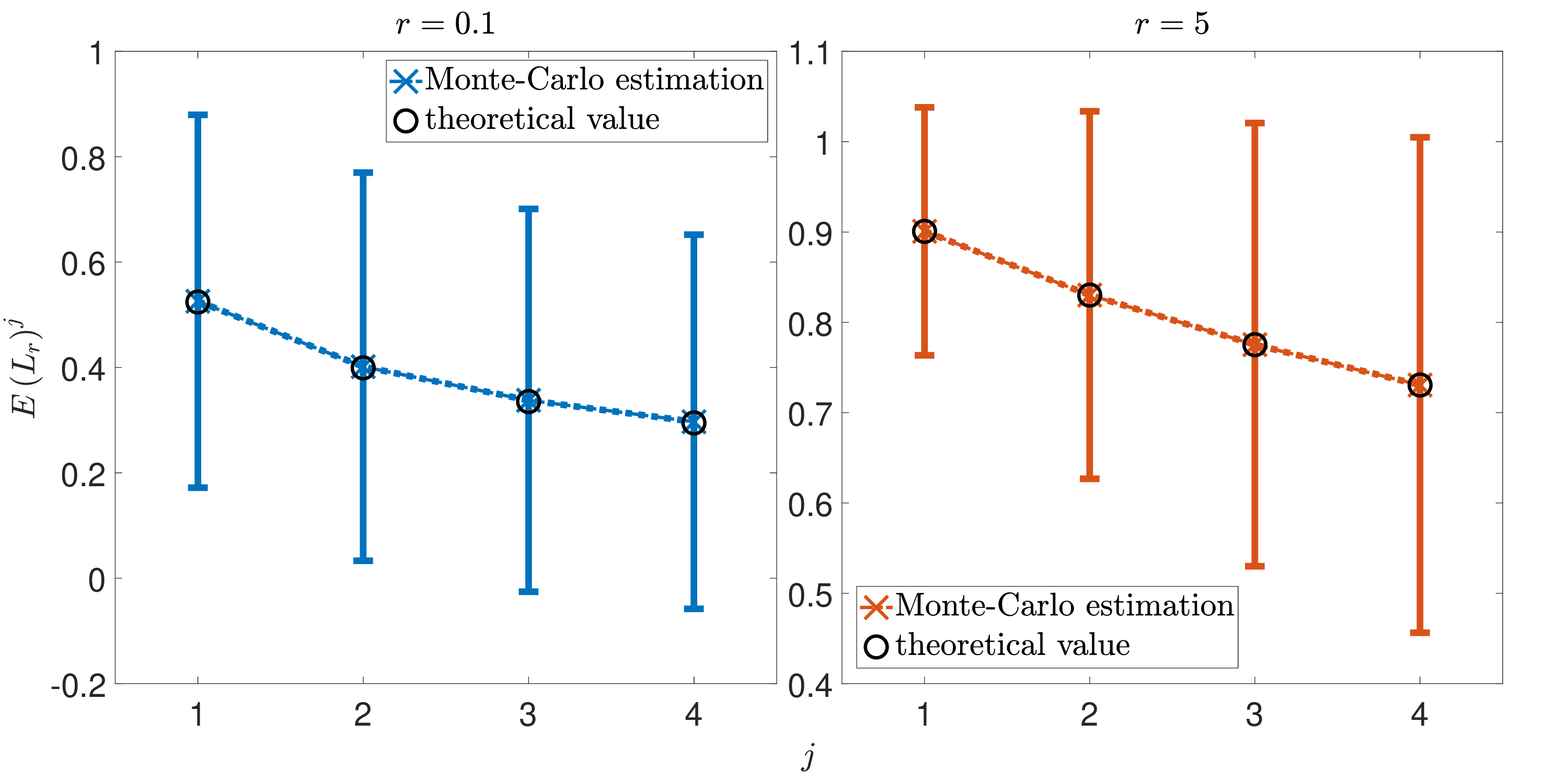}
    \caption{Monte-Carlo estimated $j$th moments of $L_r$ \eqref{eq:Lr} (crosses), their theoretical values \eqref{eq:Lmom} (circles) and the standard deviation (whiskers) of the Monte Carlo simulations.}
    \label{fig:Lmom}
\end{figure}

\begin{figure*}
    \centering
    \includegraphics[width=\linewidth]{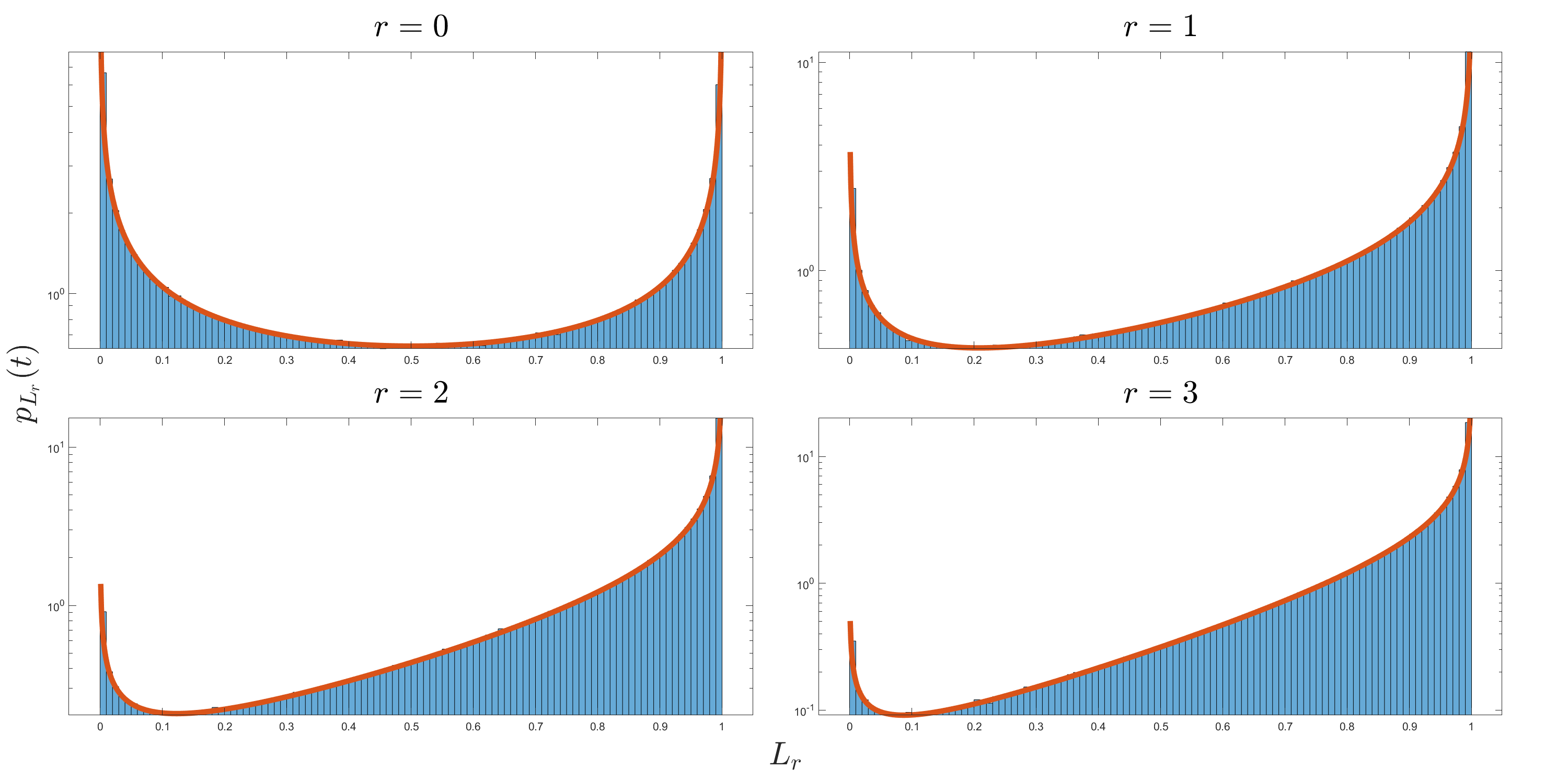}
    \caption{Monte-Carlo estimated distribution of $L_r$ \eqref{eq:Lr}, i.e. the last zero of the Wiener process under resetting with rate $r$, using Wiener process trajectories with sample size $10^5$ and timestep $10^{-4}$ compared to the analytical result \eqref{eq:LrPDF} on a logarithmic scale. Again, we observe a very good fit.}
    \label{fig:Lr}
\end{figure*}

In this section we will inspect the distribution of the last moment on the time interval $[0,1]$ in which the Wiener process under resetting hits zero. In the classical setting, for the Brownian motion $W(t)$, it is defined as
\begin{equation}
    L=\sup_t{\{t\in[0,1] : W(t)=0\}}.
\end{equation}
The second arcsine law states that
\begin{equation}
    \PRo{L\leq t}=\frac{2}{\pi}\arcsin{\left( \sqrt{t} \right)}.
\end{equation}
As in the first law, we have the scaling property
\begin{equation}
    \sup_t{\{t\in[0,\theta] : W(t)=0\}}\stackrel{D}{=}\theta \sup_t{\{t\in[0,1] : W(t)=0\}}.
\end{equation}
As before, $L_r$ will be the last time the Wiener process under Poissonian resetting with intensity $r$ hits zero, i.e.
\begin{equation}\label{eq:Lr}
    L_r=\sup_t{\{t\in[0,1] : X(t)=0\}}.
\end{equation}
Recall that $X(t)$ is the Wiener process under Poissonian resetting.
In our case the resetting point is $0$, so $L_r$ will be the time of the last zero of the Wiener process in the last interresetting interval before $t=1$. We thus have
\begin{equation}
L_r\stackrel{D}{=}\tau'_{N(1)}+(1-\tau'_{N(1)})L,
\end{equation}
where $\tau'_{N(1)}$ is the time of the last reset before $t=1$. Now using the law of total probability with respect to $N(1)$ we get that the PDF of $L_r$ equals
\begin{equation}\label{eq:Ltotal}
    p_{L_r}(t)\dd t=\sum_{k=0}^\infty\PRo{L_r\in(t,t+\dd t)|N(1)=k}\PRo{N(1)=k}.
\end{equation}
The time of the last reset before $t=1$ given $k$ resets can be easily calculated from \eqref{eq:poissonuniform}, as we have the relation \cite{ross2007}
\begin{equation}
    \tau'_k\stackrel{D}{=}\max\{U_1,\dots U_k\},
\end{equation}
where $U_i$ are iid $\mathcal{U}(0,1)$-distributed random variables. The PDF of $\tau'_k$ thus equals
\begin{equation}
    p_{\tau'_k}(\tau)=k\tau^{k-1}    
\end{equation}
for $\tau\in[0,1]$. Now let's average the term under the sum \eqref{eq:Ltotal} over possible $\tau'_k$ values to get
\begin{multline}\label{eq:Lavtau}
    \PRo{L_r\in(t,t+\dd t)|N(1)=k}=\int_0^1 p_{\tau'_k}(\tau)p_{\tau+(1-\tau)L}(t)\dd \tau\dd t\\
    =k\int_0^1 \tau^{k-1}p_{\tau+(1-\tau)L}(t)\dd \tau\dd t.
\end{multline}
Denote by $S_1$ the sum \eqref{eq:Ltotal} without the first component corresponding to $k=0$. When $k=0$ then we have no resets, and that occurs with probability $e^{-r}$. We have
\begin{align}\label{eq:LS1}
    \nonumber S_1&=\sum_{k=1}^\infty k\int_0^1 \tau^{k-1}p_{\tau+(1-\tau)L}(t)\dd \tau\dd t e^{-r}\frac{r^k}{k!}\\
    \nonumber &=re^{-r}\int_0^1p_{\tau+(1-\tau)L}(t)\sum_{k=1}^\infty\frac{(r\tau)^{k-1}}{(k-1)!}\dd \tau\dd t\\
    &=re^{-r}\int_0^1p_{\tau+(1-\tau)L}(t)e^{r\tau}\dd \tau\dd t.
\end{align}
The variable $\tau+(1-\tau)L$ is generalized arcsine distributed supported on $[\tau,1]$. Its PDF equals
\begin{equation}
    p_{\tau+(1-\tau)L}(t)=\frac{1}{\pi\sqrt{(t-\tau)(1-t)}}.
\end{equation}
Plugging this into formula \eqref{eq:LS1} for the sum $S_1$ we get
\begin{align}
    \nonumber S_1&=\frac{re^{-r}}{\pi\sqrt{1-t}}\int_0^t\frac{e^{r\tau}}{\sqrt{t-\tau}}\dd \tau\dd t=\frac{\sqrt{r}e^{r(t-1)}}{\pi\sqrt{1-t}}\int_{0}^{rt}\frac{e^{-u}}{u^{\sfrac{1}{2}}}\dd u\dd t\\
    &=\frac{\sqrt{r}e^{r(t-1)}}{\pi\sqrt{1-t}}\gamma\left(\sfrac{1}{2},rt\right)\dd t,
\end{align}
where $\gamma$ is the lower incomplete gamma function defined as
\begin{equation}
\gamma(s,x)=\int\limits_0^xt^{s-1}e^{-x}\dd x.
\end{equation}
Combining the results for $k>0$ and $k=0$ in \eqref{eq:Ltotal} we get the final result for the PDF of $L_r$
\begin{equation}\label{eq:LrPDF}
    p_{L_r}(t)=e^{-r}p_L(t)+\frac{\sqrt{r}e^{r(t-1)}}{\pi\sqrt{1-t}}\gamma\left(\sfrac{1}{2},rt\right).
\end{equation}
Comparison of Monte-Carlo simulations with the theoretical result for several chosen intensities $r$ is presented in figure \ref{fig:Lr}. The calculations are confirmed by the simulations. Skewness of the distribution for $r>0$ is caused by the resetting mechanism and is in line with intuition, as the last reset is in fact also a zero of the process. The mean of $L_r$ calculated using the PDF equals
\begin{equation}
    \E{L_r}=1+\frac{e^{-r}-1}{2r},
\end{equation}
which converges to $\sfrac{1}{2}$ when $r\to0$ and to $1$ when $r\to\infty$. The variance is
\begin{equation}
    \E{L_r^2}-\left(\E{L_r}\right)^2=\frac{2-e^{-r}(1+3r)-e^{-2r}}{4r^2},
\end{equation}
which goes to $0$ when $r\to\infty$. This means that indeed $L_r\stackrel{p}{\to}1$ when $r$ goes to infinity. Using the derived PDF we can also calculate $n$th moment. It is equal to
\begin{multline}\label{eq:Lmom}
    \E{L_r^n}=e^{-r}\Bigg( \frac{\Gamma{\left(\sfrac{1}{2}+n  \right)}}{\sqrt{\pi}\Gamma{(1+n)}}+r\Gamma{\left( \sfrac{3}{2}+n \right)}\\
    \times{}_2\bar{F}_2\left(1,\sfrac{3}{2}+n;\sfrac{3}{2} ,2+n;r\right)\Bigg).
\end{multline}
We can see the confirmation of the results using Monte-Carlo method in figure \ref{fig:Lmom}.
Observing the PDF \eqref{eq:LrPDF} we see, that it is a convex mixture of the arcsine distribution PDF and the function
\begin{equation}\label{eq:XPDF}
    p_X(r;t)=\frac{\sqrt{r}}{\pi \left( e^{r}-1 \right)}\frac{e^{rt}}{\sqrt{1-t}}\gamma{\left( \sfrac{1}{2},rt\right)}
\end{equation}
with weights $e^{-r}$ and $1-e^{-r}$ accordingly. When $r$ goes to $0$, then \eqref{eq:XPDF} converges to the PDF of Beta $\mathcal{B}(\sfrac{1}{2},\sfrac{3}{2})$ distribution
\begin{equation}
    p_X(0;t)=\frac{1}{B(\sfrac{1}{2},\sfrac{3}{2})}t^{\sfrac{1}{2}-1}(1-t)^{\sfrac{3}{2}-1}=\frac{2}{\pi}\sqrt{\frac{t}{1-t}},
\end{equation}
while for $r\to\infty$ we get $p_X(r;t)\to\delta(t-1)$.
\section{Third arcsine law}
Let us define $M$ as the moment the Wiener process reaches its maximum on $[0,1]$, that is
\begin{equation}
    M=\argmax_{t\in[0,1]}W(t).
\end{equation}
The arcsine law for the Wiener process states, that $M$ is arcsine distributed. Let's now define a random variable $M_r$, which is the moment in which Wiener process under Poissonian resetting reaches its maximum on $[0,1]$, i.e. 
\begin{equation}\label{eq:Max}
    M_r=\argmax_{t\in[0,1]}X(t).
\end{equation}
Conditioning on the event that there were $k$ resets, we have that
\begin{gather}
    \nonumber(I|N_1=k)\stackrel{D}{=}\argmax_{i\in\{1,\dots,k+1\}}Y^i_{\Delta\tau_i'},\\
    (M_r|N_1=k)\stackrel{D}{=}\tau'_{I-1}+M^I_{\Delta\tau'_I},
\end{gather}
where $Y^i_{\Delta\tau_i}$ is the maximal value of the $i$th Wiener process on the $i$th interresetting interval and $M^i_{\Delta\tau_i}$ is the location of the corresponding maximum. $I$ is the interresetting interval in which we observed a global maximum. We observe a strong correlation between the variables above, as the length of the intervals affects both of them. This strong dependence is the reason why the analytical approach from the previous sections does not yield feasible results in this case.

Finding the distribution of the random variable $M_r$ is a complex task and an open problem that we were unable to solve yet.
Therefore, we will turn to numerical analysis of the distribution of $M_r$. 

\begin{figure}
    \centering
    \includegraphics[width=\linewidth]{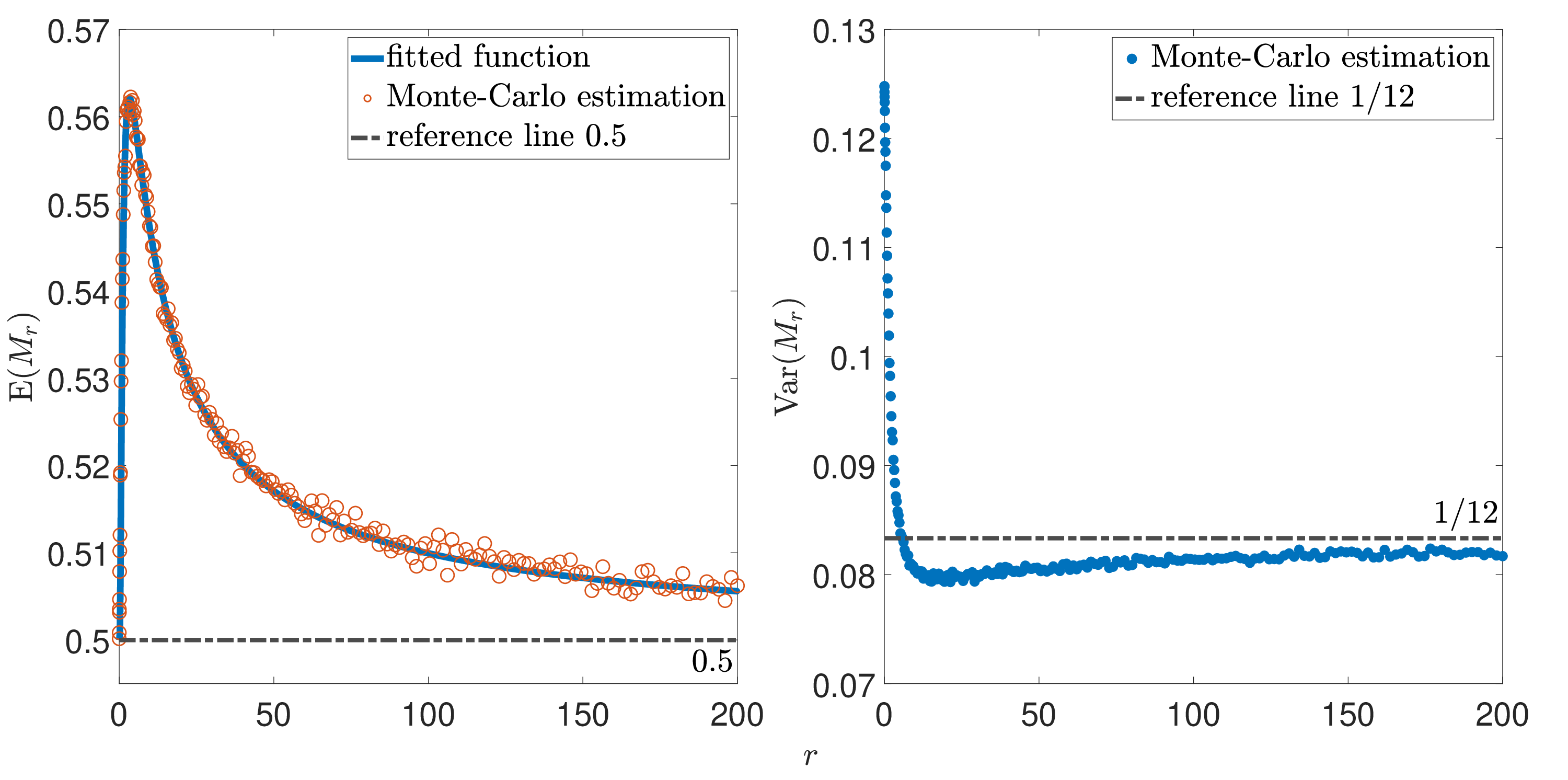}
    \caption{Mean and variance of $M_r$ \eqref{eq:Max} obtained using Monte-Carlo methods as a function of resetting rate $r$ compared with a fitted function \eqref{eq:EMrFit}. We see initial rise of the mean to the value of approximately $0.562$ at $r\approx 3.5$ before a $\sfrac{1}{r}$-like drop to the value $\sfrac{1}{2}$ being the mean of the uniform distribution $\mathcal{U}(0,1)$. Variance also drops before slow rise to the value of $\sfrac{1}{12}$.}
    \label{fig:EMr}
\end{figure}

It turns out that the mean of $M_r$ is very well fitted by the following function
\begin{equation}\label{eq:EMrFit}
    f(r)=\frac{1}{2}+\frac{ae^{-\frac{b}{r^c}}}{r^d}.
\end{equation} 
It was chosen empirically as we noticed a quick rise and steady decrease in the mean as well as expected the limiting value of $\sfrac{1}{2}$ being the mean of uniform distribution. In figure \ref{fig:EMr} we can see the mean of $M_r$ as a function of resetting rate $r$ as well as a fitted function $f(r)$ with parameters estimated using nonlinear least squares method $(a,b,c,d)=(3.3435,4.3575,0.4172,1.1309)$. The mean and variance converge to the corresponding values for a uniform random variable. We can observe that the mean attains its maximum (approximately $0.562$) at $r\approx 3.5$.
    
\begin{figure*}
    \centering    
\includegraphics[width=\linewidth]{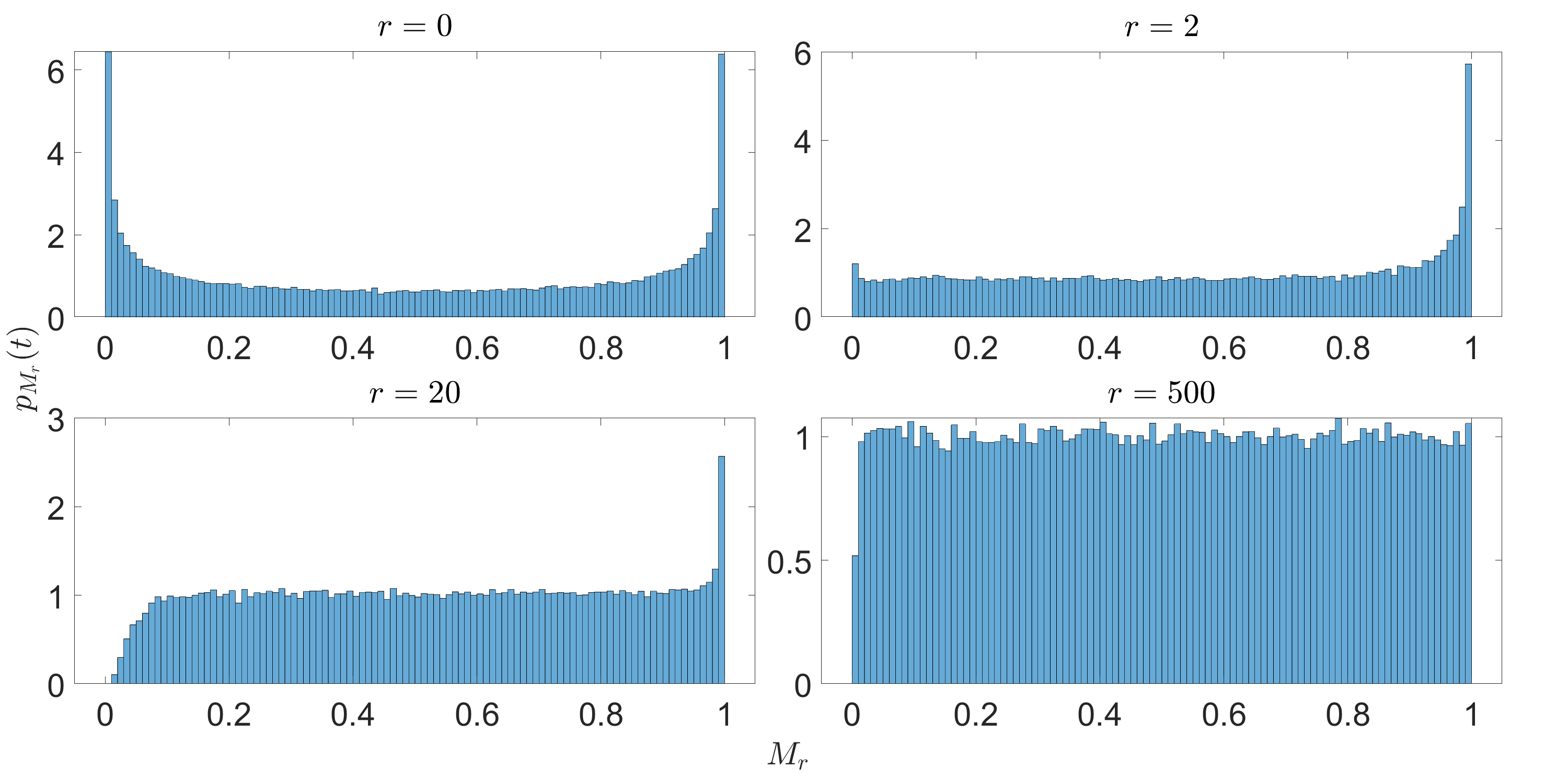}
    \caption{Histograms of simulated $M_r$ \eqref{eq:Max} variables, i.e. the times the Wiener process under resetting with rate $r$ reaches its maximum. We can see the convergence of the distribution to the uniform $\mathcal{U}(0,1)$ distribution with increasing $r$.}
    \label{fig:PMr}

\end{figure*}

In figure \ref{fig:PMr} we can see the histograms of the simulated realizations of $M_r$. The results strongly suggest the convergence of the distribution of $M_r$ to a uniform random variable with increasing $r$.
\section{Final remarks}
The question of finding a sufficiently simple expression for the distribution or the moments of the third law variable $M_r$ \eqref{eq:Max} is of great interest and provides a much greater challenge. We consider the potential follow-up research very important. The results can be generalized and used to study widely applicable jump-diffusion models.
\vfill
\bibliography{aipsamp}

\begin{thebibliography}{34}%
\makeatletter
\providecommand \@ifxundefined [1]{%
 \@ifx{#1\undefined}
}%
\providecommand \@ifnum [1]{%
 \ifnum #1\expandafter \@firstoftwo
 \else \expandafter \@secondoftwo
 \fi
}%
\providecommand \@ifx [1]{%
 \ifx #1\expandafter \@firstoftwo
 \else \expandafter \@secondoftwo
 \fi
}%
\providecommand \natexlab [1]{#1}%
\providecommand \enquote  [1]{``#1''}%
\providecommand \bibnamefont  [1]{#1}%
\providecommand \bibfnamefont [1]{#1}%
\providecommand \citenamefont [1]{#1}%
\providecommand \href@noop [0]{\@secondoftwo}%
\providecommand \href [0]{\begingroup \@sanitize@url \@href}%
\providecommand \@href[1]{\@@startlink{#1}\@@href}%
\providecommand \@@href[1]{\endgroup#1\@@endlink}%
\providecommand \@sanitize@url [0]{\catcode `\\12\catcode `\$12\catcode `\&12\catcode `\#12\catcode `\^12\catcode `\_12\catcode `\%12\relax}%
\providecommand \@@startlink[1]{}%
\providecommand \@@endlink[0]{}%
\providecommand \url  [0]{\begingroup\@sanitize@url \@url }%
\providecommand \@url [1]{\endgroup\@href {#1}{\urlprefix }}%
\providecommand \urlprefix  [0]{URL }%
\providecommand \Eprint [0]{\href }%
\providecommand \doibase [0]{http://dx.doi.org/}%
\providecommand \selectlanguage [0]{\@gobble}%
\providecommand \bibinfo  [0]{\@secondoftwo}%
\providecommand \bibfield  [0]{\@secondoftwo}%
\providecommand \translation [1]{[#1]}%
\providecommand \BibitemOpen [0]{}%
\providecommand \bibitemStop [0]{}%
\providecommand \bibitemNoStop [0]{.\EOS\space}%
\providecommand \EOS [0]{\spacefactor3000\relax}%
\providecommand \BibitemShut  [1]{\csname bibitem#1\endcsname}%
\let\auto@bib@innerbib\@empty
\bibitem [{\citenamefont {Einstein}(1905)}]{Einstein1905}%
  \BibitemOpen
  \bibfield  {author} {\bibinfo {author} {\bibfnamefont {A.}~\bibnamefont {Einstein}},\ }\href@noop {} {\bibfield  {journal} {\bibinfo  {journal} {Annalen der Physik}\ }\textbf {\bibinfo {volume} {17}},\ \bibinfo {pages} {549} (\bibinfo {year} {1905})}\BibitemShut {NoStop}%
\bibitem [{\citenamefont {Wiener}(1923)}]{Wiener1923}%
  \BibitemOpen
  \bibfield  {author} {\bibinfo {author} {\bibfnamefont {N.}~\bibnamefont {Wiener}},\ }\href@noop {} {\bibfield  {journal} {\bibinfo  {journal} {Journal of Mathematics and Physics}\ }\textbf {\bibinfo {volume} {2}},\ \bibinfo {pages} {131} (\bibinfo {year} {1923})}\BibitemShut {NoStop}%
\bibitem [{\citenamefont {Mörters}\ and\ \citenamefont {Peres}(2010)}]{Morters2010}%
  \BibitemOpen
  \bibfield  {author} {\bibinfo {author} {\bibfnamefont {P.}~\bibnamefont {Mörters}}\ and\ \bibinfo {author} {\bibfnamefont {Y.}~\bibnamefont {Peres}},\ }\href@noop {} {\emph {\bibinfo {title} {Brownian Motion}}}\ (\bibinfo  {publisher} {Cambridge University Press},\ \bibinfo {year} {2010})\BibitemShut {NoStop}%
\bibitem [{\citenamefont {Yen}\ and\ \citenamefont {Yor}(2013)}]{Yen2013}%
  \BibitemOpen
  \bibfield  {author} {\bibinfo {author} {\bibfnamefont {J.-Y.}\ \bibnamefont {Yen}}\ and\ \bibinfo {author} {\bibfnamefont {M.}~\bibnamefont {Yor}},\ }\enquote {\bibinfo {title} {Paul l{\'e}vy's arcsine laws},}\ in\ \href {\doibase 10.1007/978-3-319-01270-4_4} {\emph {\bibinfo {booktitle} {Local Times and Excursion Theory for Brownian Motion: A Tale of Wiener and It{\^o} Measures}}}\ (\bibinfo  {publisher} {Springer International Publishing},\ \bibinfo {address} {Cham},\ \bibinfo {year} {2013})\ pp.\ \bibinfo {pages} {43--54}\BibitemShut {NoStop}%
\bibitem [{\citenamefont {Alsmeyer}\ and\ \citenamefont {Buckmann}(2019)}]{Alsmeyer2019}%
  \BibitemOpen
  \bibfield  {author} {\bibinfo {author} {\bibfnamefont {G.}~\bibnamefont {Alsmeyer}}\ and\ \bibinfo {author} {\bibfnamefont {F.}~\bibnamefont {Buckmann}},\ }\href {\doibase https://doi.org/10.1016/j.spa.2018.02.014} {\bibfield  {journal} {\bibinfo  {journal} {Stochastic Processes and their Applications}\ }\textbf {\bibinfo {volume} {129}},\ \bibinfo {pages} {223} (\bibinfo {year} {2019})}\BibitemShut {NoStop}%
\bibitem [{\citenamefont {Kotani}(2006)}]{Kotani2006}%
  \BibitemOpen
  \bibfield  {author} {\bibinfo {author} {\bibfnamefont {S.}~\bibnamefont {Kotani}},\ }\href@noop {} {\bibfield  {journal} {\bibinfo  {journal} {Lecture Notes in Mathematics}\ }\textbf {\bibinfo {volume} {1874}},\ \bibinfo {pages} {333} (\bibinfo {year} {2006})}\BibitemShut {NoStop}%
\bibitem [{\citenamefont {Yor}(1992)}]{Yor1992}%
  \BibitemOpen
  \bibfield  {author} {\bibinfo {author} {\bibfnamefont {M.}~\bibnamefont {Yor}},\ }\href {https://books.google.pl/books?id=d8OFwQEACAAJ} {\emph {\bibinfo {title} {Some Aspects of Brownian Motion Part I : Some Special Functionals}}},\ Lectures in Mathematics. ETH Z{\"u}rich\ (\bibinfo  {publisher} {Birkh{\"a}user Basel},\ \bibinfo {year} {1992})\BibitemShut {NoStop}%
\bibitem [{\citenamefont {Spitzer}(1956)}]{Spitzer1956}%
  \BibitemOpen
  \bibfield  {author} {\bibinfo {author} {\bibfnamefont {F.}~\bibnamefont {Spitzer}},\ }\href@noop {} {\bibfield  {journal} {\bibinfo  {journal} {Transactions of the American Mathematical Society}\ }\textbf {\bibinfo {volume} {82}},\ \bibinfo {pages} {323} (\bibinfo {year} {1956})}\BibitemShut {NoStop}%
\bibitem [{\citenamefont {Knight}(1963)}]{Knight1963}%
  \BibitemOpen
  \bibfield  {author} {\bibinfo {author} {\bibfnamefont {F.~B.}\ \bibnamefont {Knight}},\ }\href@noop {} {\bibfield  {journal} {\bibinfo  {journal} {Transactions of the American Mathematical Society}\ }\textbf {\bibinfo {volume} {109}},\ \bibinfo {pages} {56} (\bibinfo {year} {1963})}\BibitemShut {NoStop}%
\bibitem [{\citenamefont {Borodin}\ and\ \citenamefont {Salminen}(2002)}]{Borodin2002}%
  \BibitemOpen
  \bibfield  {author} {\bibinfo {author} {\bibfnamefont {A.~N.}\ \bibnamefont {Borodin}}\ and\ \bibinfo {author} {\bibfnamefont {P.}~\bibnamefont {Salminen}},\ }\href@noop {} {\emph {\bibinfo {title} {Handbook of Brownian Motion - Facts and Formulae}}}\ (\bibinfo  {publisher} {Birkhäuser},\ \bibinfo {year} {2002})\BibitemShut {NoStop}%
\bibitem [{\citenamefont {Gallardo}\ and\ \citenamefont {Pérez-Castillo}(2021)}]{Gallardo2021}%
  \BibitemOpen
  \bibfield  {author} {\bibinfo {author} {\bibfnamefont {D.}~\bibnamefont {Gallardo}}\ and\ \bibinfo {author} {\bibfnamefont {I.}~\bibnamefont {Pérez-Castillo}},\ }\href@noop {} {\bibfield  {journal} {\bibinfo  {journal} {Physical Review E}\ }\textbf {\bibinfo {volume} {104}},\ \bibinfo {pages} {034123} (\bibinfo {year} {2021})}\BibitemShut {NoStop}%
\bibitem [{\citenamefont {Karatzas}\ and\ \citenamefont {Shreve}(1991)}]{Karatzas1991}%
  \BibitemOpen
  \bibfield  {author} {\bibinfo {author} {\bibfnamefont {I.}~\bibnamefont {Karatzas}}\ and\ \bibinfo {author} {\bibfnamefont {S.}~\bibnamefont {Shreve}},\ }\href {\doibase https://doi.org/10.1007/978-1-4612-0949-2} {\emph {\bibinfo {title} {Brownian Motion and Stochastic Calculus}}},\ Graduate Texts in Mathematics (113)\ (\bibinfo  {publisher} {Springer New York},\ \bibinfo {year} {1991})\BibitemShut {NoStop}%
\bibitem [{\citenamefont {Ding}\ and\ \citenamefont {Wu}(2019)}]{Ding2019}%
  \BibitemOpen
  \bibfield  {author} {\bibinfo {author} {\bibfnamefont {J.}~\bibnamefont {Ding}}\ and\ \bibinfo {author} {\bibfnamefont {Y.}~\bibnamefont {Wu}},\ }\href@noop {} {\bibfield  {journal} {\bibinfo  {journal} {Journal of Theoretical Probability}\ }\textbf {\bibinfo {volume} {32}},\ \bibinfo {pages} {2025} (\bibinfo {year} {2019})}\BibitemShut {NoStop}%
\bibitem [{\citenamefont {Revuz}\ and\ \citenamefont {Yor}(1999)}]{Revuz1999}%
  \BibitemOpen
  \bibfield  {author} {\bibinfo {author} {\bibfnamefont {D.}~\bibnamefont {Revuz}}\ and\ \bibinfo {author} {\bibfnamefont {M.}~\bibnamefont {Yor}},\ }\href@noop {} {\emph {\bibinfo {title} {Continuous Martingales and Brownian Motion}}}\ (\bibinfo  {publisher} {Springer},\ \bibinfo {year} {1999})\BibitemShut {NoStop}%
\bibitem [{\citenamefont {Billingsley}(1995)}]{Billingsley1995}%
  \BibitemOpen
  \bibfield  {author} {\bibinfo {author} {\bibfnamefont {P.}~\bibnamefont {Billingsley}},\ }\href@noop {} {\emph {\bibinfo {title} {Probability and Measure}}}\ (\bibinfo  {publisher} {Wiley},\ \bibinfo {year} {1995})\BibitemShut {NoStop}%
\bibitem [{\citenamefont {Breen}\ and\ \citenamefont {West}(2022)}]{Breen2022}%
  \BibitemOpen
  \bibfield  {author} {\bibinfo {author} {\bibfnamefont {K.}~\bibnamefont {Breen}}\ and\ \bibinfo {author} {\bibfnamefont {B.~J.}\ \bibnamefont {West}},\ }\href@noop {} {\bibfield  {journal} {\bibinfo  {journal} {Physical Review E}\ }\textbf {\bibinfo {volume} {105}},\ \bibinfo {pages} {044123} (\bibinfo {year} {2022})}\BibitemShut {NoStop}%
\bibitem [{\citenamefont {Evans}\ and\ \citenamefont {Majumdar}(2011)}]{Evans2011}%
  \BibitemOpen
  \bibfield  {author} {\bibinfo {author} {\bibfnamefont {M.~R.}\ \bibnamefont {Evans}}\ and\ \bibinfo {author} {\bibfnamefont {S.~N.}\ \bibnamefont {Majumdar}},\ }\href@noop {} {\bibfield  {journal} {\bibinfo  {journal} {Physical Review Letters}\ }\textbf {\bibinfo {volume} {106}},\ \bibinfo {pages} {160601} (\bibinfo {year} {2011})}\BibitemShut {NoStop}%
\bibitem [{\citenamefont {Evans}\ \emph {et~al.}(2020)\citenamefont {Evans}, \citenamefont {Majumdar},\ and\ \citenamefont {Schehr}}]{Evans2020}%
  \BibitemOpen
  \bibfield  {author} {\bibinfo {author} {\bibfnamefont {M.~R.}\ \bibnamefont {Evans}}, \bibinfo {author} {\bibfnamefont {S.~N.}\ \bibnamefont {Majumdar}}, \ and\ \bibinfo {author} {\bibfnamefont {G.}~\bibnamefont {Schehr}},\ }\href {\doibase 10.1088/1751-8121/ab7cfe} {\bibfield  {journal} {\bibinfo  {journal} {Journal of Physics A: Mathematical and Theoretical}\ }\textbf {\bibinfo {volume} {53}},\ \bibinfo {pages} {193001} (\bibinfo {year} {2020})}\BibitemShut {NoStop}%
\bibitem [{\citenamefont {A.~Pal}\ and\ \citenamefont {Evans}(2015)}]{Pal2015}%
  \BibitemOpen
  \bibfield  {author} {\bibinfo {author} {\bibfnamefont {A.~K.}\ \bibnamefont {A.~Pal}}\ and\ \bibinfo {author} {\bibfnamefont {M.~R.}\ \bibnamefont {Evans}},\ }\href@noop {} {\bibfield  {journal} {\bibinfo  {journal} {Journal of Physics A: Mathematical and Theoretical}\ }\textbf {\bibinfo {volume} {49}},\ \bibinfo {pages} {225001} (\bibinfo {year} {2015})}\BibitemShut {NoStop}%
\bibitem [{\citenamefont {Reuveni}(2016)}]{Reuveni2016}%
  \BibitemOpen
  \bibfield  {author} {\bibinfo {author} {\bibfnamefont {S.}~\bibnamefont {Reuveni}},\ }\href@noop {} {\bibfield  {journal} {\bibinfo  {journal} {Physical Review Letters}\ }\textbf {\bibinfo {volume} {116}},\ \bibinfo {pages} {170601} (\bibinfo {year} {2016})}\BibitemShut {NoStop}%
\bibitem [{\citenamefont {Saikat~Ray}\ and\ \citenamefont {Reuveni}(2020)}]{Ray2020}%
  \BibitemOpen
  \bibfield  {author} {\bibinfo {author} {\bibfnamefont {D.~M.}\ \bibnamefont {Saikat~Ray}}\ and\ \bibinfo {author} {\bibfnamefont {S.}~\bibnamefont {Reuveni}},\ }\href@noop {} {\bibfield  {journal} {\bibinfo  {journal} {Journal of Physical Chemistry Letters}\ }\textbf {\bibinfo {volume} {11}},\ \bibinfo {pages} {7769} (\bibinfo {year} {2020})}\BibitemShut {NoStop}%
\bibitem [{\citenamefont {Pal}\ and\ \citenamefont {Prasad}(2017)}]{Pal2017}%
  \BibitemOpen
  \bibfield  {author} {\bibinfo {author} {\bibfnamefont {A.}~\bibnamefont {Pal}}\ and\ \bibinfo {author} {\bibfnamefont {V.~V.}\ \bibnamefont {Prasad}},\ }\href@noop {} {\bibfield  {journal} {\bibinfo  {journal} {Physical Review E}\ }\textbf {\bibinfo {volume} {96}},\ \bibinfo {pages} {062135} (\bibinfo {year} {2017})}\BibitemShut {NoStop}%
\bibitem [{\citenamefont {Benjamin~Besga}\ and\ \citenamefont {Majumdar}(2020)}]{Besga2020}%
  \BibitemOpen
  \bibfield  {author} {\bibinfo {author} {\bibfnamefont {A.~B.}\ \bibnamefont {Benjamin~Besga}}\ and\ \bibinfo {author} {\bibfnamefont {S.~N.}\ \bibnamefont {Majumdar}},\ }\href@noop {} {\bibfield  {journal} {\bibinfo  {journal} {Physical Review Letters}\ }\textbf {\bibinfo {volume} {124}},\ \bibinfo {pages} {200601} (\bibinfo {year} {2020})}\BibitemShut {NoStop}%
\bibitem [{\citenamefont {Łukasz Kuśmierz}\ \emph {et~al.}(2014)\citenamefont {Łukasz Kuśmierz}, \citenamefont {Majumdar},\ and\ \citenamefont {Redner}}]{Kusmierz2014}%
  \BibitemOpen
  \bibfield  {author} {\bibinfo {author} {\bibnamefont {Łukasz Kuśmierz}}, \bibinfo {author} {\bibfnamefont {S.~N.}\ \bibnamefont {Majumdar}}, \ and\ \bibinfo {author} {\bibfnamefont {S.}~\bibnamefont {Redner}},\ }\href@noop {} {\bibfield  {journal} {\bibinfo  {journal} {Physical Review Letters}\ }\textbf {\bibinfo {volume} {113}},\ \bibinfo {pages} {220602} (\bibinfo {year} {2014})}\BibitemShut {NoStop}%
\bibitem [{\citenamefont {Kundu}\ and\ \citenamefont {Reuveni}(2024)}]{Kundu2024}%
  \BibitemOpen
  \bibfield  {author} {\bibinfo {author} {\bibfnamefont {A.}~\bibnamefont {Kundu}}\ and\ \bibinfo {author} {\bibfnamefont {S.}~\bibnamefont {Reuveni}},\ }\href {\doibase 10.1088/1751-8121/ad1e1b} {\bibfield  {journal} {\bibinfo  {journal} {Journal of Physics A: Mathematical and Theoretical}\ }\textbf {\bibinfo {volume} {57}},\ \bibinfo {pages} {060301} (\bibinfo {year} {2024})}\BibitemShut {NoStop}%
\bibitem [{\citenamefont {Eliazar}\ and\ \citenamefont {Reuveni}(2020)}]{Eli2020}%
  \BibitemOpen
  \bibfield  {author} {\bibinfo {author} {\bibfnamefont {I.}~\bibnamefont {Eliazar}}\ and\ \bibinfo {author} {\bibfnamefont {S.}~\bibnamefont {Reuveni}},\ }\href {\doibase 10.1088/1751-8121/abae8c} {\bibfield  {journal} {\bibinfo  {journal} {Journal of Physics A: Mathematical and Theoretical}\ }\textbf {\bibinfo {volume} {53}} (\bibinfo {year} {2020}),\ 10.1088/1751-8121/abae8c}\BibitemShut {NoStop}%
\bibitem [{\citenamefont {Eliazar}\ and\ \citenamefont {Reuveni}(2021{\natexlab{a}})}]{Eli2021}%
  \BibitemOpen
  \bibfield  {author} {\bibinfo {author} {\bibfnamefont {I.}~\bibnamefont {Eliazar}}\ and\ \bibinfo {author} {\bibfnamefont {S.}~\bibnamefont {Reuveni}},\ }\href {\doibase 10.1088/1751-8121/ac16c5} {\bibfield  {journal} {\bibinfo  {journal} {Journal of Physics A: Mathematical and Theoretical}\ }\textbf {\bibinfo {volume} {54}},\ \bibinfo {pages} {355001} (\bibinfo {year} {2021}{\natexlab{a}})}\BibitemShut {NoStop}%
\bibitem [{\citenamefont {Eliazar}\ and\ \citenamefont {Reuveni}(2021{\natexlab{b}})}]{Eli2021_2}%
  \BibitemOpen
  \bibfield  {author} {\bibinfo {author} {\bibfnamefont {I.}~\bibnamefont {Eliazar}}\ and\ \bibinfo {author} {\bibfnamefont {S.}~\bibnamefont {Reuveni}},\ }\href {\doibase 10.1088/1751-8121/abe4a0} {\bibfield  {journal} {\bibinfo  {journal} {Journal of Physics A: Mathematical and Theoretical}\ }\textbf {\bibinfo {volume} {54}} (\bibinfo {year} {2021}{\natexlab{b}}),\ 10.1088/1751-8121/abe4a0}\BibitemShut {NoStop}%
\bibitem [{\citenamefont {Eliazar}\ and\ \citenamefont {Reuveni}(2023{\natexlab{a}})}]{Eli2023}%
  \BibitemOpen
  \bibfield  {author} {\bibinfo {author} {\bibfnamefont {I.}~\bibnamefont {Eliazar}}\ and\ \bibinfo {author} {\bibfnamefont {S.}~\bibnamefont {Reuveni}},\ }\href {\doibase 10.1088/1751-8121/acb183} {\bibfield  {journal} {\bibinfo  {journal} {Journal of Physics A: Mathematical and Theoretical}\ }\textbf {\bibinfo {volume} {56}},\ \bibinfo {pages} {024002} (\bibinfo {year} {2023}{\natexlab{a}})}\BibitemShut {NoStop}%
\bibitem [{\citenamefont {Eliazar}\ and\ \citenamefont {Reuveni}(2023{\natexlab{b}})}]{Eli2023_2}%
  \BibitemOpen
  \bibfield  {author} {\bibinfo {author} {\bibfnamefont {I.}~\bibnamefont {Eliazar}}\ and\ \bibinfo {author} {\bibfnamefont {S.}~\bibnamefont {Reuveni}},\ }\href {\doibase 10.1088/1751-8121/acb184} {\bibfield  {journal} {\bibinfo  {journal} {Journal of Physics A: Mathematical and Theoretical}\ }\textbf {\bibinfo {volume} {56}},\ \bibinfo {pages} {024003} (\bibinfo {year} {2023}{\natexlab{b}})}\BibitemShut {NoStop}%
\bibitem [{\citenamefont {Ghosh}\ \emph {et~al.}(2023)\citenamefont {Ghosh}, \citenamefont {Nayak}, \citenamefont {Liu}, \citenamefont {Li},\ and\ \citenamefont {Marchesoni}}]{Ghosh2023}%
  \BibitemOpen
  \bibfield  {author} {\bibinfo {author} {\bibfnamefont {P.~K.}\ \bibnamefont {Ghosh}}, \bibinfo {author} {\bibfnamefont {S.}~\bibnamefont {Nayak}}, \bibinfo {author} {\bibfnamefont {J.}~\bibnamefont {Liu}}, \bibinfo {author} {\bibfnamefont {Y.}~\bibnamefont {Li}}, \ and\ \bibinfo {author} {\bibfnamefont {F.}~\bibnamefont {Marchesoni}},\ }\href {\doibase 10.1063/5.0159148} {\bibfield  {journal} {\bibinfo  {journal} {The Journal of Chemical Physics}\ }\textbf {\bibinfo {volume} {159}},\ \bibinfo {pages} {031101} (\bibinfo {year} {2023})},\ \Eprint {http://arxiv.org/abs/https://pubs.aip.org/aip/jcp/article-pdf/doi/10.1063/5.0159148/18050097/031101\_1\_5.0159148.pdf} {https://pubs.aip.org/aip/jcp/article-pdf/doi/10.1063/5.0159148/18050097/031101\_1\_5.0159148.pdf} \BibitemShut {NoStop}%
\bibitem [{\citenamefont {Feller}(1957)}]{feller1957}%
  \BibitemOpen
  \bibfield  {author} {\bibinfo {author} {\bibfnamefont {W.}~\bibnamefont {Feller}},\ }\href {https://books.google.pl/books?id=NFNQAAAAMAAJ} {\emph {\bibinfo {title} {An Introduction to Probability Theory and Its Applications, Volume 2}}},\ \bibinfo {series} {An Introduction to Probability Theory and Its Applications}\ No.\ \bibinfo {number} {t. 1-2}\ (\bibinfo  {publisher} {Wiley},\ \bibinfo {year} {1957})\BibitemShut {NoStop}%
\bibitem [{\citenamefont {Ross}(2007)}]{ross2007}%
  \BibitemOpen
  \bibfield  {author} {\bibinfo {author} {\bibfnamefont {S.}~\bibnamefont {Ross}},\ }\href {https://books.google.pl/books?id=-dwtxagQ9skC} {\emph {\bibinfo {title} {Introduction to Probability Models}}},\ Introduction to Probability Models\ (\bibinfo  {publisher} {Elsevier Science},\ \bibinfo {year} {2007})\BibitemShut {NoStop}%
\bibitem [{\citenamefont {Abramowitz}\ and\ \citenamefont {Stegun}(1964)}]{Abramowitz1964}%
  \BibitemOpen
  \bibfield  {author} {\bibinfo {author} {\bibfnamefont {M.}~\bibnamefont {Abramowitz}}\ and\ \bibinfo {author} {\bibfnamefont {I.}~\bibnamefont {Stegun}},\ }\href@noop {} {\emph {\bibinfo {title} {Handbook of Mathematical Functions}}}\ (\bibinfo  {publisher} {United States Department of Commerce, National Bureau of Standards},\ \bibinfo {year} {1964})\BibitemShut {NoStop}%
\end{thebibliography}%

\end{document}